\RequirePackage{fix-cm}
\documentclass[smallextended]{svjour3}

\smartqed

\usepackage{amssymb,amsmath}
\usepackage{graphicx,color}
\usepackage{url}
\usepackage{bm}
\usepackage[noend]{algpseudocode}
\usepackage{algorithm}
\usepackage[colorinlistoftodos]{todonotes}
\usepackage{multirow}

\usepackage{hyperref}
\hypersetup{
  colorlinks=true,
  urlcolor=blue,
  linkcolor=blue,
  citecolor=blue
}

\usepackage[]{natbib}
\setcitestyle{aysep={}}

\usepackage{graphicx}
\usepackage{amssymb}
\usepackage{amscd}
\usepackage{amsmath}
\usepackage[noend]{algpseudocode}
\usepackage{algorithm}
\usepackage{enumerate}
\usepackage[all]{xy}
\usepackage{comment}
\usepackage{subfigure}
\usepackage{color}

\newcommand{\R}{{\mathbb{R}}}
\newcommand{\X}{{\mathbb{X}}}
\newcommand{\Z}{{\mathbb{Z}}}

\newcommand{\cech}{{\rm \check{C}ech}}
\newcommand{\rips}{{\rm Rips}}
\journalname{Journal name here}
\begin{document}

\title{Persistence Diagrams with Linear Machine Learning Models\thanks{
    This work is partially supported by 
    JSPS KAKENHI Grant Number JP 16K17638, 
    JST CREST Mathematics15656429, 
    JST ``Materials research by Information Integration" Initiative (${\rm MI}^2{\rm I}$) project of the Support Program for Starting Up Innovation Hub, 
    Structural Materials for Innovation Strategic Innovation Promotion Program D72, and New Energy and Industrial Technology Development Organization (NEDO).
  }}


\author{Ippei Obayashi \and Yasuaki Hiraoka}

\institute{Ippei Obayashi \at
  Advanced Institute for Materials Research (WPI-AIMR), Tohoku University. \\
  2-1-1 Katahira, Aoba-ku, Sendai, 980-8577 Japan \\
  Tel.: +81-22-217-6320\\
  Fax: +81-22-217-5129\\
  \email{ippei.obayashi.d8@tohoku.ac.jp}           
  \and
  Y. Hiraoka \at
  Advanced Institute for Materials Research (WPI-AIMR), Tohoku University.\\
  Center for Materials research by Information Integration (CMI2), 
  Research and Services Division of Materials Data and Integrated System   (MaDIS), National Institute for Materials Science (NIMS). \\
\email{hiraoka@tohoku.ac.jp}           
}

\date{Received: date / Accepted: date}

\maketitle

\begin{abstract}
Persistence diagrams have been widely recognized as a compact descriptor for characterizing multiscale topological features in data. When many datasets are available, statistical features embedded in those persistence diagrams can be extracted by applying machine learnings. In particular, the ability for explicitly analyzing the inverse in the original data space from those statistical features of persistence diagrams is significantly important for practical applications. In this paper, we propose a unified method for the inverse analysis by combining linear machine learning models with persistence images. The method is applied to point clouds and cubical sets, showing the ability of the statistical inverse analysis and its advantages. 

\end{abstract}

\section{Introduction}
\label{sec:intro}
Given a dataset, its statistical features can be extracted by applying machine learning methods~\citep{prml}. 
Needless to say, machine learning is now one of the central scientific and engineering subjects, and is rapidly enlarging its theoretical foundations and ranges of practical applications. 
For example, in materials science, the amount of available data has recently been increasing due to improvement of experimental methods and computational resources. These datasets are expected to be used for further developments of high performance materials based on machine learnings, leading to a new concept called ``materials informatics"  \citep{rajan1,rajan2}.

As another branch of data science, topological data analysis (TDA)~\citep{carlsson,eh} has also been rapidly developed from theoretical aspects to applications in the last decade. In TDA, persistent homology and its persistence diagram \citep{elz,zc} are widely used  for capturing multiscale topological features in data. Recent improvements of efficient computations of persistence diagrams \citep{dipha,phat} enable us to apply them into practical problems such as materials science \citep{amorphous,granular,ichinomiya,iron}, sensor networks \citep{sensor}, evolutions of virus~\citep{virus} etc. 
As a descriptor of data, persistence diagrams have the following significant properties: translation and rotation invariance, multi-scalability, and robustness for noise. Together with developments of statistical foundations \citep{landscape,chazal,fasy,kusano1,kusano2,pssk}, persistence diagrams nowadays have been recognized as a compact descriptor for complicated data. 

In a series of works on materials TDA \citep{amorphous,granular,ichinomiya,iron}, analyzing the inverse in the original data space (atomic configurations or digital images) from persistence diagrams is significantly important to explicitly study the materials structures and properties. Therefore,  toward further progress that materials TDA incorporates with materials informatics, we need to develop a framework of machine learnings on persistence diagrams which allows the inverse analysis.

In this paper, we propose a unified method for studying the shape of data by using persistence diagrams with machine learnings in both direct and inverse problems. The essence of our method is to combine persistence images ~\citep{persistence_image} and linear machine learning models. 

For standard machine learning methods, the input data is supposed to be given by vectors, and therefore we need to transform persistence diagrams into vectors. Some vectorization methods of persistence diagrams have been proposed in the literatures~\citep{persistence_image,landscape,kusano1,kusano2,pssk}, and we here use persistence images. This is because it allows us to reconstruct persistence diagrams from vectors obtained by machine learning results, providing a key step in the inverse route of our analysis. 

Taking this advantage, we apply linear models of machine learnings to persistence images. Since the learned result of linear machine learning models is given by a vector with the same dimension as input vectors, we can easily reconstruct the persistence diagram from the learned result using the functionality of persistence images. 
Namely, the persistence diagram itself is obtained as learning. Furthermore, by studying inverse problems from the reconstructed persistence diagram to the original data space, we can explicitly characterize statistically significant topological features embedded in data. In this paper, we deal with an inverse problem studying the locations of birth-death pairs of persistence diagrams in the original data space. As another advantage using linear machine learning models, we also propose an important concept called sparse persistence diagram. This new concept allows us to discard irrelevant generators and to focus on most significant ones in the reconstructed persistence diagram for learning tasks.

It should be remarked that, for only direct problems such as predictions from data, nonlinear methods such as kernel methods and neural networks are possibly appropriate, because such nonlinear transformations often make the prediction performance better than linear models. However, if our interest is to understand mechanisms of data structures, the inverse route going back to the original data from the learned results is inevitable. 

As summary, the contribution of this paper is to propose a unified method in topological data analysis with the ability to study inverse problems by combining the following methods:
  \begin{enumerate}
  \item Persistence images
  \item Linear machine learning models
  \item Inverse analysis of persistence diagrams
  \end{enumerate}
In Section \ref{sec:method}, after brief introduction of our input data formats and persistent homology, we recall persistence images and linear models of machine learnings used in this paper. Then, in Section \ref{sec:result}, we demonstrate our method to some problems on point clouds and cubical sets, and show the performance comparing to other methods. Some future problems and related topics including some practical applications in materials science are summarized in Section \ref{sec:conclusion}.

\section{Methods}\label{sec:method}
We first explain some preliminaries about geometric models and persistent homology. Although the theory of persistent homology has been rapidly extended in various general settings, we here introduce the minimum necessary for later discussions. Readers who want to understand the theory in higher generality are encouraged to study the latest literatures. 

\subsection{Geometric models}\label{sec:geometric_model}
In this paper, we mainly consider two types of input data. The first type is given by a finite points $P=\{x_i\in\R^N\colon i=1,\dots,m\}$ in a Euclidean space $\R^N$, which is also called a \emph{point cloud} in TDA. For example, this data type is frequently used for expressing atomic configurations in materials science. 

Our interest is to characterize multiscale topological properties in $P$, and to this aim, we consider the $r$-ball model
\[
	P_r=\bigcup_{i=1}^mB_r(x_i),
\]
where $B_r(x_i)=\{y\in\R^N\colon ||y-x_i||\leq r\}$ is the ball with radius $r$ centered at $x_i$. By construction, when the radius $r$ is very small (resp. large), $P_r$ has the same topology as $N$ disconnected points (resp. one point). Between these two extremal cases, $P_r$ may exhibit appearance and disappearance of holes by changing the radius $r$. Note that we have a natural inclusion $P_r\subset P_s$ for $r\leq s$, meaning that the radius parameter $r$ can be regarded as a resolution of the point cloud $P$. 

For practical data analysis, the $r$-ball model $P_r$ is not convenient to handle in computers, and hence we usually build  simplicial complex models from $P_r$. For instance, the \emph{{\v C}ech complex} $\cech(P,r)$ and the \emph{Rips complex} (or Vietoris-Rips complex) $\rips(P,r)$ are simplicial complexes with the vertex set $P$ whose $k$-simplex is assigned by the following rule, respectively,
\begin{align*}
&\{x_{i_0},\dots,x_{i_k}\}\in \cech(P,r)\Leftrightarrow \bigcap_{s=0}^kB_r(x_{i_s})\neq \emptyset,\\
&\{x_{i_0},\dots,x_{i_k}\}\in \rips(P,r)\Leftrightarrow B_r(x_{i_s})\cap B_r(x_{i_t})\neq\emptyset,~0\leq \forall s<\forall t\leq k.
\end{align*}
Note that, by construction, both simplicial complex models naturally define a (right continuous) \emph{filtration}. Namely, for $X_r=\cech(P,r)$ or $X_r=\rips(P,r)$, it satisfies $X_r\subset X_s$ for $r\leq s$ and $X_s=\bigcap_{s<t}X_t$. In this section, we denote the filtration by $\X=\{X_r\colon r\in\R\}$.

Our next data type is given by a cubical set, which is a standard mathematical expression for digital images. Following the notation used in the reference \citep{kmm}, let $I\subset \R$ be an elementary interval, i.e., 
\[
	I=[\ell,\ell+1]\quad {\rm or}\quad I=[\ell,\ell]
\]
for some $\ell\in\Z$. An elementary cube $Q=I_1\times\dots\times I_N\subset \R^N$ is defined by a product of elementary intervals $I_i$. Then, a subset $X\subset \R^N$ is said to be cubical if $X$ can be expressed as a union of elementary cubes in $\R^N$. 

Let us denote by $\mathcal{K}^N_W$ the set of all elementary cubes in the window $\Lambda_W=[-W,W]^N\subset \R^N$. 
Given a function $f:\mathcal{K}_W^N\rightarrow \R$, we can build a cubical set in $\Lambda_W$  as a sublevel set
\begin{align}\label{ref:sublevel}
X_t=\bigcup\{Q\in\mathcal{K}^N_W\colon f(Q)\leq t\}
\end{align}
for each parameter $t$. In practical applications such as digital image analysis, this function is often given by the Manhattan distance (e.g., see Figure \ref{fig:ph_binary_image}) or a grayscale function. It is easy to see that the cubical sets $X_t$ also lead to a filtration $\X=\{X_t\colon t\in\R\}$.

These are the two standard types of our input data.  We note that those filtrations satisfy the properties that $X_t=\emptyset$ for sufficiently small $t$ and $X_t$ is acyclic\footnote{A topological space $X$ with $\tilde H_q(X)=0$ for any $q$ is called acyclic, where $\tilde H_q(X)$ is the reduced homology of $X$.} for sufficiently large $t$, respectively.

\subsection{Persistent homology}
Let $\Bbbk$ be a field. In this paper, the $q$th homology $H_q(X)$ of a topological space $X$ is defined over the field $\Bbbk$, and hence $H_q(X)$ is given as a $\Bbbk$-vector space. Intuitively, the dimension of $H_q(X)$ as a $\Bbbk$-vector space counts the number of $q$-dimensional holes in $X$, and each basis vector expresses the corresponding $q$-dimensional hole in $X$, where, for example, $q=0,1,2$ express connected components, rings, and cavities, respectively. 
Then, given a pair of topological spaces $X\hookrightarrow Y$, we can define the induced linear map $\varphi: H_q(X)\rightarrow H_q(Y)$, which characterizes whether a hole in $X$ persists in $Y$ or not. 

The input to the persistent homology is given by a filtration $\X=\{X_t\colon t\in \R\}$ of topological spaces. 
In this paper, $X_t$ is given by a simplicial complex or a cubical set.
For simplicity, we also assume the properties for filtrations remarked in the final paragraph in Section \ref{sec:geometric_model}, although we do not really need it by modifying the argument here. 
Then, the $q$th \emph{persistent homology} $H_q(\X)=(H_q(X_t),\varphi_s^t)$ of the filtration $\X$ is defined by the family of homologies $\{H_q(X_t)\colon t\in\R\}$ and the induced linear maps $\varphi_s^t: H_q(X_s)\rightarrow H_q(X_t)$ for all $s\leq t$. 

Under the assumption of our filtrations, the persistent homology $H_q(\X)$ can be uniquely decomposed by using the so-called interval representations:
\begin{align}\label{eq:decomposition}
H_q(\X)\simeq \bigoplus_{i=1}^pI(b_i,d_i),
\end{align}
where $b_i,d_i\in\R$ with $b_i<d_i$.
Here, $I(b_i,d_i)=(U_t,f_s^t)$ consists of a family of vector spaces
\begin{align*}
U_t=\left\{\begin{array}{ll}
\Bbbk,&b_i\leq t<d_i,\\
0,&{\rm otherwise,}
\end{array}\right.
\end{align*}
and the identity map $f_s^t={\rm id}_\Bbbk$ for $b_i\leq s\leq t < d_i$. Note that the $0$th persistent homology in (\ref{eq:decomposition}) is understood as the reduced sense, meaning that one connected component which persists for any large $t\in\R$ is removed. Each interval representation $I(b_i,d_i)$ is also  called a generator of $H_q(\X)$. 

Each generator $I(b_i,d_i)$ expresses that a $q$-dimensional hole appears in $\X$ at the parameter $t=b_i$, persists up to $t<d_i$, and then disappears at $t=d_i$. We call $b_i, d_i, d_i-b_i$ the birth time, death time, and lifetime of $I(b_i,d_i)$, respectively. 

Under the unique decomposition (\ref{eq:decomposition}), the $q$th persistence diagram $D_q(\X)$ of $\X$ is defined by a multiset\footnote{A multiset is a set with multiplicity of each point.}
\begin{align*}
D_q(\X)=\{(b_i,d_i)\in \Delta\colon i=1,\dots,p\},
\end{align*}
where $\Delta=\{(b,d)\in\R^2\colon b < d\}$. It is known that the birth-death pair $(b_i,d_i)\in D_q(\X)$ with large lifetime can be regarded as reliable topological structure in $\X$, while that with small lifetime is likely to be a noisy structure. This statement is justified by the stability theorem of persistent homology \citep{ceh}.

For a review about computational aspect of persistent homology, we refer the readers to the paper \citep{review}.

\subsection{Examples}
Here, we show several examples to make clear the concepts explained so far.  To this aim, the examples are chosen to be simple enough for demonstration.

\begin{figure}[htbp]
  \centering
  \includegraphics[width=\hsize]{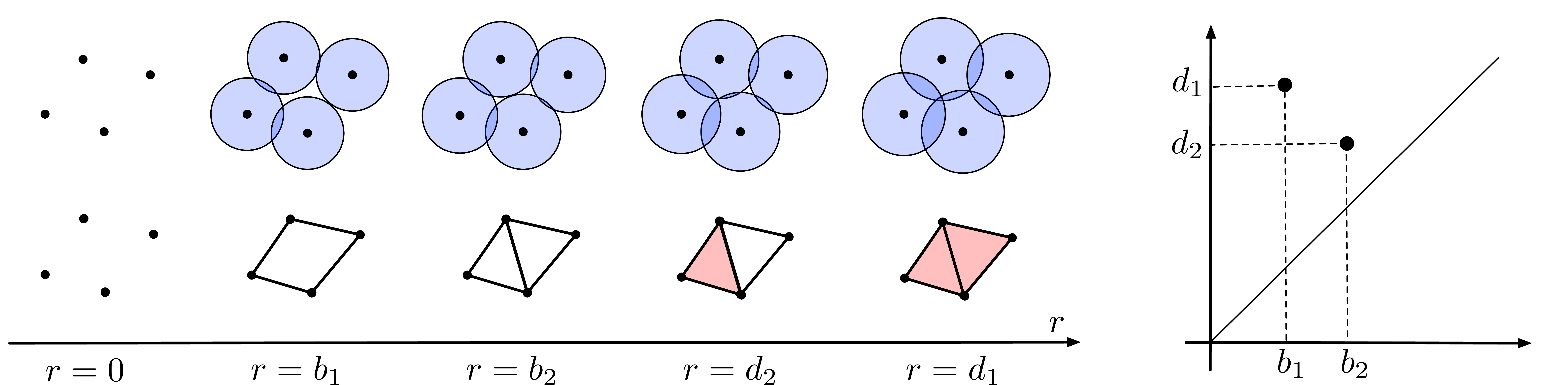}
  \caption{Left top: Filtration of the $r$-ball models. Left bottom: Filtration of the corresponding {\v C}ech complexes. Right: The $1$st persistence diagram.
  (1) A ring is bone at $r=b_1$. 
  (2) Another ring is bone at $r=b_2$. 
  (3) The second ring dies at $r=d_2$.  
  (4) The first ring dies at $r=d_1$.
  }
  \label{fig:four_points}
\end{figure}

We first consider an example of a point cloud given by four points on the plane shown in the left ($r=0$) of Figure~\ref{fig:four_points}. As we explained, each point is replaced by a ball and we study topological changes during the fattening process of the balls by increasing the radii. This fattening process is drawn on the left top of Figure~\ref{fig:four_points}, while the sequence below expresses its {\v C}ech complex filtration. 

At the radius $r=b_1$, the first ring is born, and we record its birth parameter as $b_1$. Similarly, the second ring appears at the birth parameter  $r=b_2$. On the other hand, at radius $r=d_1,d_2$, those rings disappear and we record them as their death parameters. Hence, the $1$st persistence diagram of the {\v C}ech complex filtration is given by $\{(b_1, d_1), (b_2, d_2)\}$, which is shown on the right of Figure~\ref{fig:four_points}. 

\begin{figure}[htbp]
  \centering
  \includegraphics[width=\hsize]{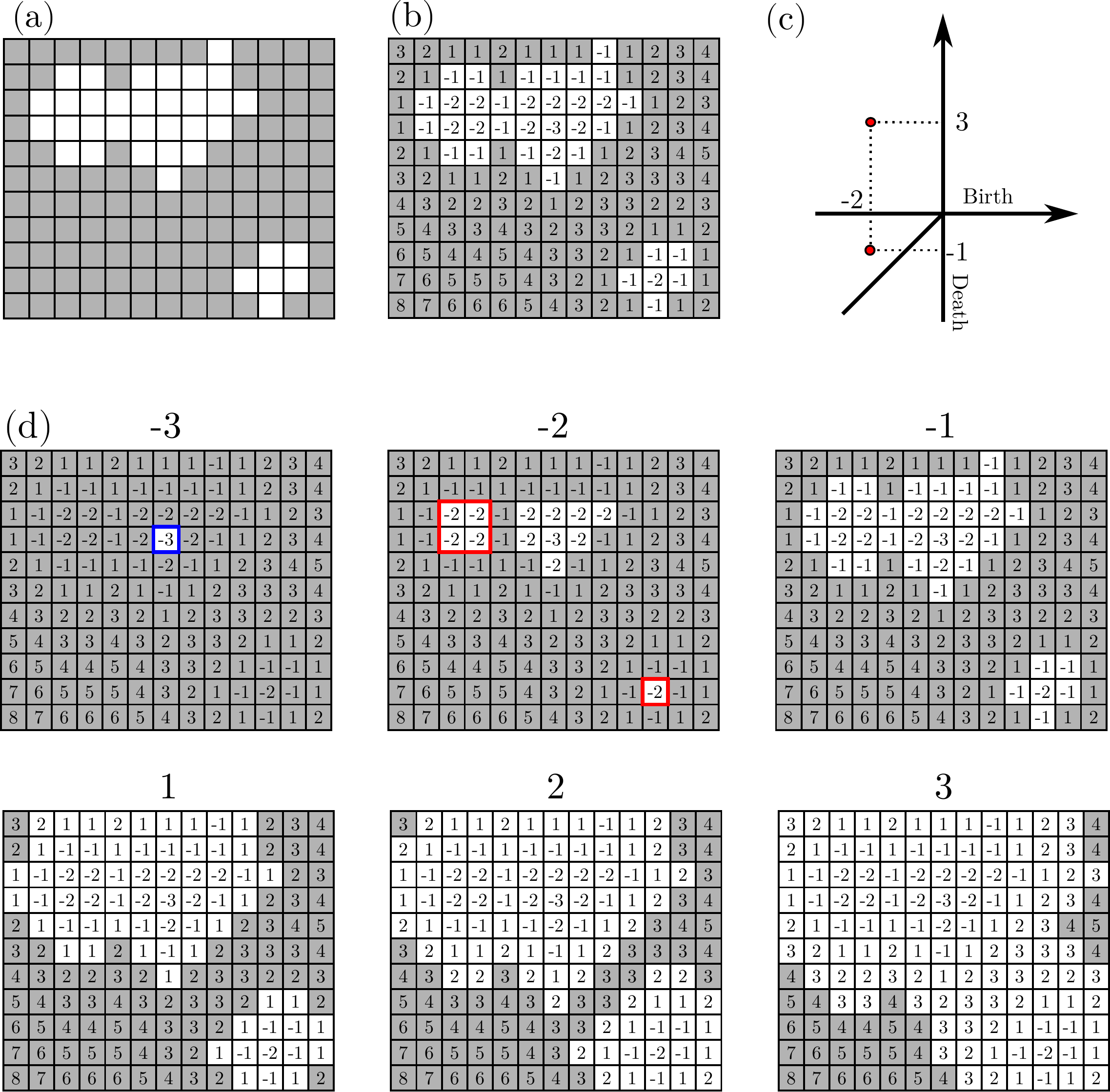}
  \caption{
  (a) Input binary image. (b) Manhattan distance. (c) $0$th reduced persistence diagram. (d) Filtration of binary images with respect to the Manhattan distance. 
  The colored squares indicate the initial locations of three connected components.
The blue square indicates the connected component removed in the reduced persistent homology.}
  \label{fig:ph_binary_image}
\end{figure}

Next, we consider an example of a cubical set. The input data is given by a binary image (a) in Figure~\ref{fig:ph_binary_image}, and we consider a function $f$ assigning an integer on each pixel shown in (b). 
Here, positive (resp. negative) numbers are assigned to the gray (resp. white) pixels using the Manhattan distance. Then, following the construction (\ref{ref:sublevel}) of sublevel sets, we obtain a filtration of cubical sets of white pixels shown in (d). 

\begin{figure}[htbp]
  \centering
  \includegraphics[width=0.5\hsize]{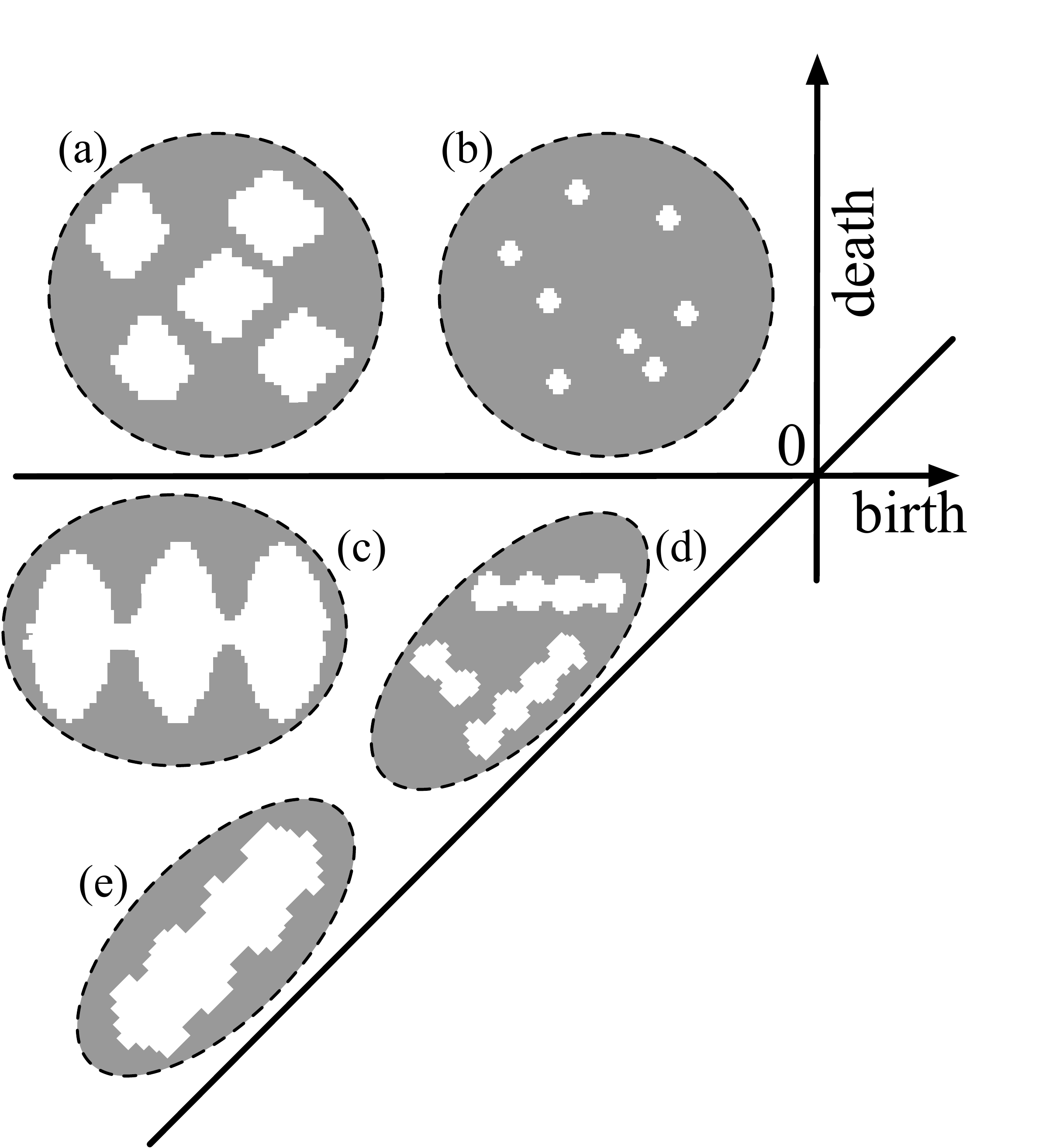}  
  \caption{Cubical sets drawn in the dashed circles (a)-(e) express typical geometric structures in the $0$th persistence diagram. (a) Large islands. (b) Small islands. (c) Large islands with narrow bridges. (d) Narrow bands. (e) Broad bands. For (a) and (b), the births correspond to the radii of the islands. For (d) and (e), the births and deaths correspond to the half widths of the bands. }
  \label{fig:pd_binary_image_cheatsheet}
\end{figure}

In this example, three connected components appear in the filtration and those birth events are colored in blue and red. The death of those generators corresponds to a merging to another connected component. Then, the 0th reduced persistence diagram is given by $\{(-2, -1), (-2, 3)\}$,  shown in (c). Note that the first connected component born at $-3$ is removed in the reduced persistence diagram. 
We also note that, from the assignment $f$ using Manhattan distance, all birth parameters take negative values. Figure \ref{fig:pd_binary_image_cheatsheet} summarizes common geometric structures captured by the 0th persistence diagram based on the Manhattan distance.

For deep analysis using persistence diagrams, we often want to know the origin of each birth-death pair. One easy and useful way is to utilize a death simplex (resp. death cube) for a point cloud (resp. cubical set).  In the {\v C}ech filtration model of Figure~\ref{fig:four_points}, two generators (i.e., rings) die when each red simplex fills the corresponding ring, and those simplices show the locations of the generators. We call these locations \emph{death positions} of the generators. Even for generators with higher dimensions and also for the setting of cubical sets, this idea works in a similar way. 

On the other hand, for generators with dimension zero, birth events may possess the information of locations. In Figure \ref{fig:ph_binary_image}(d), there are two red squares and they express the central locations of each connected component. We call these locations \emph{birth positions} of the corresponding birth-death pairs. 

We note that the birth/death positions are easily obtained in standard algorithms of computing persistence diagrams, and hence no additional computations are required. These techniques, which will be demonstrated in the later section, are exploited for some practical analysis in materials science \citep{iron}. It should be remarked that, if one wants to obtain further information about the inverse of birth-death pairs, the technique of optimal cycles \citep{optimal1,optimal2} can be another choice, although it requires a much computer resource.

\subsection{Persistence images}
Recall that a persistence diagram is a multi-set on $\R^2$. Hence, we need to vectorize persistence diagrams to apply  machine learning models. In this paper, we use the \emph{persistence image}~\citep{persistence_image} for vectorization. 

Given a $q$th persistence diagram $D_{q} = \{(b_k,d_k)\in\Delta\colon k=1,\dots,\ell\}$, the persistence image $\rho$ is defined by a function on $\R^2$ as
\begin{align}
  \rho(x, y) &= \sum_{k=1}^\ell w(b_k, d_k) \exp\left(-\frac{(b_k-x)^2 + (d_k - y)^2}{2 \sigma^2}\right), \nonumber\\
  w(b, d) &= \arctan ( C (d - b)^p). \label{eq:weight}
\end{align}
Here, $C > 0$, $p > 0$, $\sigma > 0$ are parameters, $w(b, d)$ is a 
weight function, and we regard the function $\rho$ as a vector in a function space $L^2(\R^2)$. We remark that the weight function is chosen so that we can respect the significance of generators according to its lifetimes in the statistical analysis. As we see in Section \ref{sec:result}, the parameters are usually determined to be appropriate values using cross validations. 

For computations, we discretize the persistence image $\rho$ and construct a histogram on the plane with an appropriate finite mesh. 
Obviously, since all birth-death pairs are located in 
$\{(b, d) \mid N_- \leq b < d \leq N_+ \}$ with some constants $N_-,N_+$, the histogram is constructed on this area.
Then, we obtain a vector from the discretization of $\rho$ by ordering the elements on the grids in a prefixed order. Note that the dimension of the  vector is equal to the number of grids used for the histogram. In the following, we also call the discretization of $\rho$ the persistence image.

We note that there are several methods for vectorizations of persistence diagrams. One important advantage using persistence images is that we can easily reconstruct a histogram from a vector, and hence can obtain a corresponding persistence diagram. 
However, it is not straightforward in general to reconstruct persistence diagrams from vectors in nonlinear vectorizations.  This advantage is effectively used in our method. 

We also remark that, precisely speaking, the weight function (\ref{eq:weight}) is not used in the original paper \citep{persistence_image} but first studied in the paper \citep{kusano2}, in which  performance comparisons with different weights for persistence images and also with other vectorizations are thoroughly discussed. For details, we refer the readers to the paper \citep{kusano2}.

\subsection{Linear machine learning models}\label{sec:linear_models}
In this section, we briefly recall the logistic regression and the linear regression as standard supervised machine learning methods \citep{prml}.

In the linear regression model,
we consider a pair of an input vector $x \in \R^n$ (called \emph{explanatory variable})
and its output value $y \in \R$ (called \emph{response variable}), and study the relation between them in the linear form
\begin{align*}
  y = w \cdot x + b + \textrm{(noise)},
\end{align*}
where $w \in \R^n$ and $b \in \R$ are unknown parameters and the
noise is randomly determined from a normal distribution with mean 0.
From a set of known input-output pairs $\{(x_i, y_i)\}_{i=1}^M$, called
a \emph{training set}, we find an optimal $w$ and $b$ for the model.
Such optimal parameters are derived by minimizing the following \emph{mean squared loss error function} with respect to $w$ and $b$:
\begin{align}
  E(w, b) = \frac{1}{2M} \sum_{i=1}^M (w\cdot x_i+b - y_i)^2 .\label{eq:error-linear-regression}
\end{align}

In the logistic regression model for a binary classification task,
we consider a pair of an input vector $x \in \R^n$ and its output value $y\in\{0,1\}$, and study the relation of classification $0/1$ based on the following form
\begin{equation}
  \label{eq:logistic}
  \begin{aligned}
    P(y = 1 \mid w, b) &= g(w\cdot x + b), \\
    P(y = 0 \mid w, b) &= 1 - P(y=1 \mid w, b) = g(-w\cdot x - b), \\
    g(z) &= 1/(e^{-z}+1),
  \end{aligned}
\end{equation}
where $w \in \R^n$ and $b \in \R$ are unknown parameters.
From training data $\{(x_i, y_i)\}_{i=1}^M$, we find an optimal $w$ and
$b$ in a similar way to the linear regression. Here, optimal parameters are given
by minimizing the following \emph{cross entropy error function}:
\begin{equation}
  \label{eq:error-logistic-regression}
  \begin{aligned}
    L(w, b) &= -\frac{1}{M}
      \sum_{i=1}^M\left\{y_i \log \hat{y}_i + (1-y_i) \log (1 - \hat{y}_i) \right\},\\
    \hat{y}_i &= g(w \cdot x_i + b).
  \end{aligned}
\end{equation}

We note that, for both the linear regression and logistic regression, these optimization problems are equivalent to the maximization of the log likelihood. 

In our method, the input vector $x \in \R^n$ is given by vectorized persistence diagrams using persistence images. Then, the learned vector $w$ becomes a dual vector to the persistence images, and especially, its dimension  is the same as $x$. Hence, $w$ can be expressed as a (dual) persistence diagram by the reverse process of the vectorization using persistence images. In this way, our method outputs persistence diagrams as learning results. 

For practical applications, we often encounter the problem of
over-fittings, 
if the dimension $n$ of input vectors is relatively large compared to the data size $M$. 
Under this condition, the result of the optimization problem excessively fit the training set and does not give appropriate performance for untrained data. 
In our setting, since the dimension of vectors obtained from persistence images is very large, we usually face the over-fitting problem.
The vectors given by
persistence images also have another statistical problem called multicollinearity\citep{statistics}.
Adjacent grids elements of a vector by persistence image
are strongly correlated because of Gaussian diffusion, and such a strong correlation
causes the difficulty of determining coefficients and the numerical instability.

One effective way for avoiding the over-fitting and multicollinearity is to add a regularization (penalty) term $R(w)$ into the error function.
Namely, we minimize the following modified error functions for $w$ and $b$\begin{align*}
  E(w, b) + \lambda R(w) &\ \ \mbox{(for a linear regression)},\\
  L(w, b) + \lambda R(w) &\ \ \mbox{(for a logistic regression),}
\end{align*}
where $\lambda>0$ is a weight parameter controlling the regularization effect. 
Typical regularization terms are given as
\begin{align*}
  R(w) = \frac{1}{2} \|w\|^2_2,   \quad R(w) = \|w\|_1. \\
\end{align*}
The former is called an $\ell^2$-regularization and the latter is called an
$\ell^1$-regularization. A linear regression with the $\ell^2$-regularization is called ridge, while a linear regression with the $\ell^1$-regularization
is called lasso~\citep{lasso}.

The advantage of the $\ell^2$-regularization is its good mathematical property.
For example, $\ell^2$-regularization term is differentiable but $\ell^1$-regularization term is not.
The ridge optimization problem has the closed form solution. However, the lasso does not have such forms.

On the other hand, the $\ell^1$-regularization has a significant property of the \emph{sparsity}. A vector $w$ is called sparse if its elements are all zero except for only a few elements. It is well-known that the learned vector $w$ under $\ell^1$-regularization becomes a sparse vector, and hence we obtain a \emph{sparse persistence diagram} as a  result of learning. 
As we will see later, the sparseness of the learned persistence diagram is often very useful, when we interpret the learned results.

The parameter $\lambda$ of the regularization term controls the complexity of the learned result~\citep{prml}. When the weight $\lambda$ becomes larger, the regularization term $R(w)$ becomes smaller. This means that $w$ becomes more sparse in the $\ell^1$-regularization. Such a reduction of the complexity is useful for finding the most essential elements for regressions. However, 
when $\lambda$ is too large, the learned results may drop important information. Therefore, we need to determine a suitable $\lambda$ in practice. A validation set or cross validation method are often applied to choose such a parameter~\citep{prml}.  The effect of changing $\lambda$ in our method is discussed in Section~\ref{sec:result}.

\subsection{Summary of our method}

\begin{enumerate}
\item Prepare an input data $\{(g_i, y_i)\}_{i=1}^M$. Here, each $g_i$ is a point cloud or a digital image, and $y_i$ is a real value for the linear regression or 0/1 value for the logistic regression.
\item Compute the persistence diagram $D^{(i)}$ from $g_i$.
\item Compute the vectorization $x_i \in \R^n$ of $D^{(i)}$ using the persistence image. 
\item Apply the linear regression or the logistic regression with a regularization term to the data $\{(x_i, y_i)\}_{i=1}^M$ and find $w \in \R^n$ and $b \in \R$.
Choose the $\ell^2$- or $\ell^1$-regularization, depending on the purpose. 
\item The learned result $w$ is visualized by the reconstruction of the persistence diagram from $w$.
From the reconstructed persistence diagram, one may extract important areas on the diagram. 
\item For explicitly identifying the geometric structure of those important areas on the diagram, one can study the birth/death positions. 
\end{enumerate}

\section{Results and discussions}\label{sec:result}
In this section, we demonstrate the performance of our methods for logistic regressions and linear regressions with binary images and point clouds. 
Here, we use filtrations of cubical sets using Manhattan distance (black: positive, white: negative) for binary images, while {\v C}ech  complex filtrations are applied for point clouds. All examples are experimented using scikit-learn\footnote{\url{http://scikit-learn.org/}} and  HomCloud\footnote{\url{http://www.wpi-aimr.tohoku.ac.jp/hiraoka_labo/homcloud-english.html}}.

\subsection{Logistic regression on binary images - an easy example}\label{sec:easy}
\begin{figure}[htbp]
  \centering
  \includegraphics[width=\hsize]{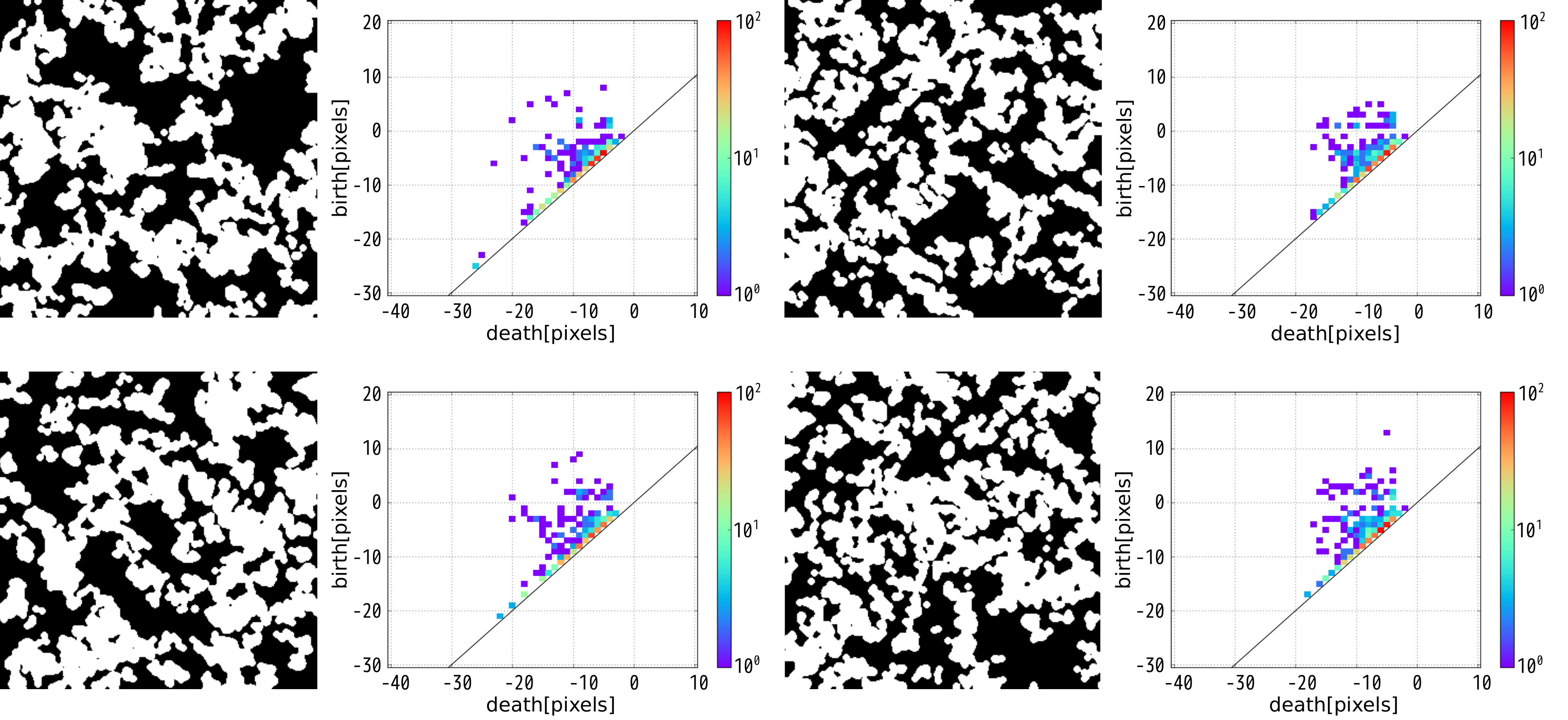}
  \caption{Input binary images and their 0th persistence diagrams.
    The left and right two images are sampled from the parameter pairs (A) and (B), respectively. }
  \label{fig:input_im_A}
\end{figure}

First, we examine the logistic regression on persistence diagrams of binary images. 
Here, the binary image data is randomly generated by 
Algorithm~\ref{alg:randomimage} in Appendix~\ref{sec:randomimage}, where 
a pair of parameters $(N,S)$ is used to generate two types of images.
One pair (A) is set to be $N=100, S=30$ and the other (B) is $N=250, S=10$. Figure~\ref{fig:input_im_A} shows the samples of both data (left: (A), right: (B)). We may intuitively observe that the images from (B) have somewhat finer structures than the images from (A). Our task is the classification of the parameters (A) and (B) from images, where we assign $0$ and $1$ for (A) and (B), respectively. 

For each parameter pair, 300 images are generated (total 2$\times$300) and 
200 of these images are sampled as a training set (total 2$\times$200). 
Then, 2$\times$100 remaining images are used as a test set to evaluate the learned result. 
Here, 0th persistence diagrams are applied for the task. The parameters of the persistence images are set to be $\sigma=2.0, C=0.5, p = 1.0$ and the mesh for the discretized persistence images is obtained by dividing the rectangle $[-40.5, 10.5]\times[-30.5,20.5]$ into $51\times 51$ grids. The $\ell^2$-regularization is used and the weight parameter $\lambda$ of the regularization term is determined by the cross validation.

In this example, the score, evaluated as the mean accuracy, of the learned result is 1.0, that is, we can perfectly identify the parameter pairs behind the images. In fact, we could also distinguish  these two parameter pairs by simply counting the number of connected components, if we had this prior knowledge. 
In Section~\ref{sec:logreg_im_difficult}, we examine a more sophisticated classification problem. For a while, let us use this example in order to explain some properties of our method.

\begin{figure}[htbp]
  \centering
  \includegraphics[width=\hsize]{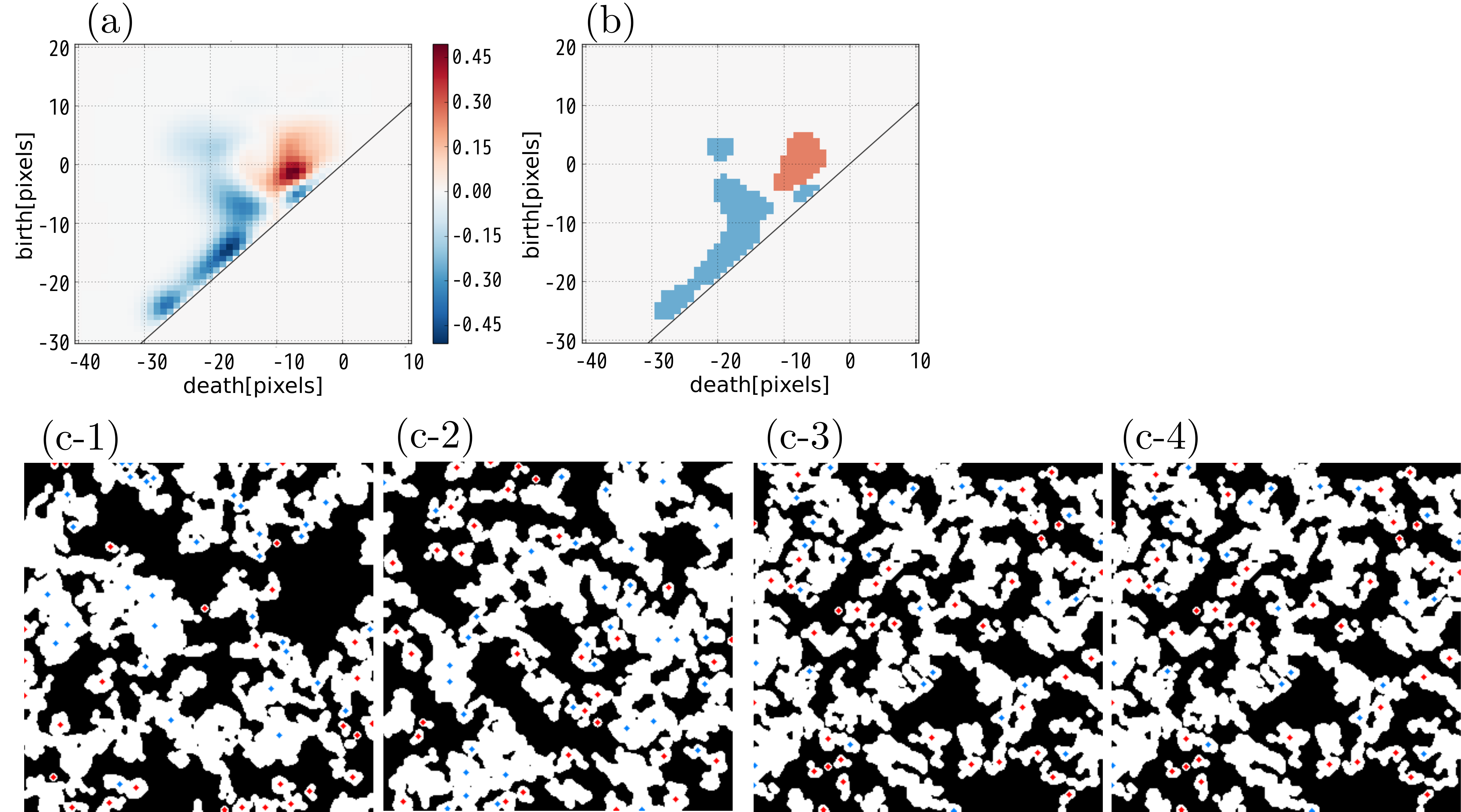}
  \caption{(a) The reconstructed persistence diagram from the learned vector $w$. The blue (resp. red) area contributes to the class 0 (resp. 1). (b) A thresholding of (a). (c-1)$\sim$(c-4) The birth positions 
  of the generators in blue and red areas in (b) are plotted with the same color. }
  \label{fig:output_im_A_1}
\end{figure}

Figure~\ref{fig:output_im_A_1} (a) shows the reconstructed persistence diagram from the learned vector $w$, and (b) shows the area at which the magnitude is above a certain threshold. Recalling the classification rule (\ref{eq:logistic}), nonzero elements in $w$ (and hence nonzero generators in its reconstructed persistence diagram) work for making classification decisions. Namely, from the 0/1 assignment rule, the reconstructed persistence diagram concludes that generators in the blue (resp. red) area statistically contribute to the classification (A) (resp. (B)).

Furthermore, by plotting the birth positions of these generators, we can explicitly identify the geometric structures which characterize the parameter pairs. Figure (c-1)$\sim$(c-4) show those birth positions, where the blue (resp. red) points correspond to the blue (resp. red) area. Recalling the interpretation in Figure \ref{fig:pd_binary_image_cheatsheet}, we find that the characteristic geometric structures of (B) are explained by small islands and narrow bands  whose inner radii are $4\sim 10$ pixels; this is consistent to our intuition that (B) contains finer structures.

\begin{figure}[htbp]
  \centering
  \includegraphics[width=\hsize]{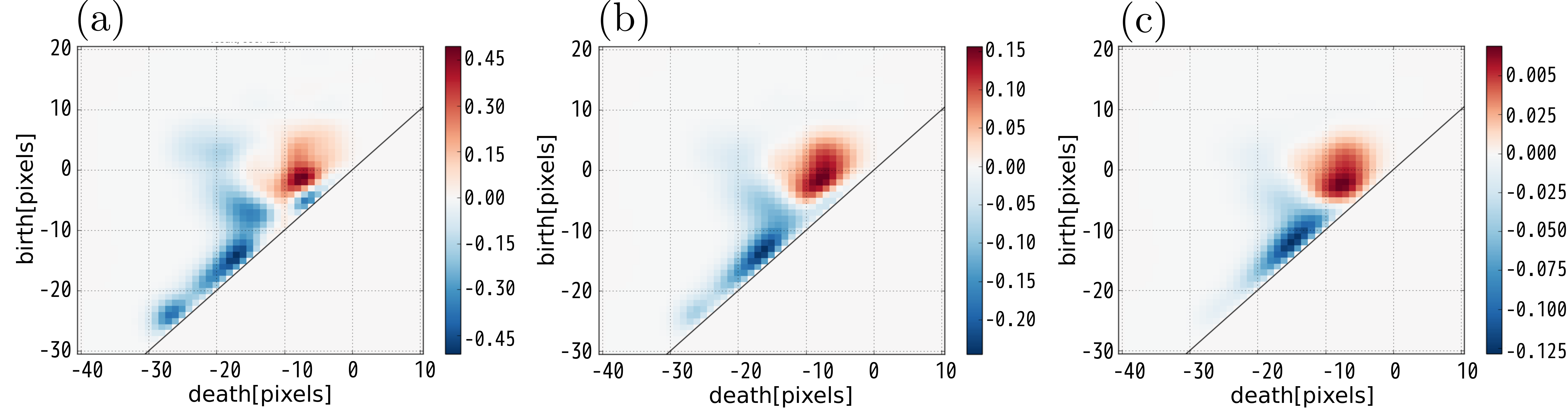}
  \caption{The reconstructed persistence diagrams with several weight parameters $\lambda$.
    (a) $\lambda = 0.35938$ (determined by the cross validation), (b) $\lambda=10$, and  (c) $\lambda = 100$.
  }
  \label{fig:output_im_A_l2}
\end{figure}

Using this example, let us study the effect of the weight parameter $\lambda$ for the regularization. Figure~\ref{fig:output_im_A_l2} shows the reconstructed persistence diagrams from the learned vectors for several weight parameters $\lambda$. When $\lambda$ becomes larger, in addition to the fact that the magnitude of the persistence diagram becomes smaller, its distribution becomes simpler. This is because the weigh parameter $\lambda$ of the regularization controls the complexity of the learned result, which is expressed in the distribution of the reconstructed persistence diagram.

\begin{figure}[htbp]
  \centering
  \includegraphics[width=\hsize]{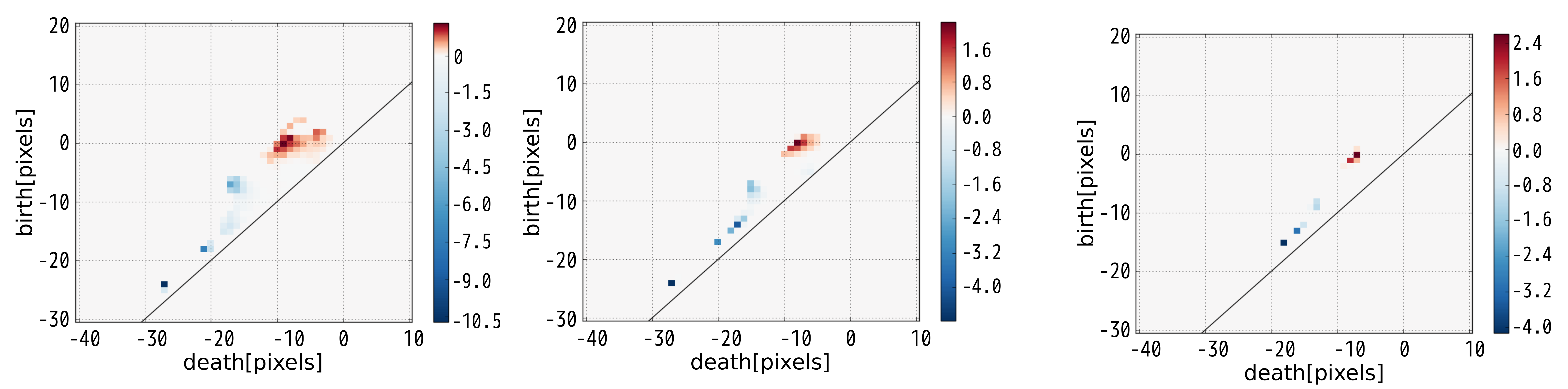}
  \caption{The reconstructed persistence diagrams using the $\ell^1$-regularization with several weight parameters $\lambda$.
  (a) $\lambda=0.01$, (b) $\lambda = 0.1$, (c) $\lambda = 1$.}
  \label{fig:output_im_A_l1}
\end{figure}

We also compare with the $\ell^1$-regularization in this example. 
Figure~\ref{fig:output_im_A_l1} shows the reconstructed persistence diagrams  using the $\ell^1$-regularization with several parameters $\lambda$.
As mentioned in Section \ref{sec:linear_models}, an important property of the $\ell^1$-regularization is the sparseness of the learned result $w$. In our method, this property is reflected as sparse persistence diagram. Hence, again recalling the classification rule (\ref{eq:logistic}), the selected few grids in the sparse persistence diagram are supposed to work most effectively for the classification task. In other words, the birth-death pairs around the grids are especially important for the classification.
Furthermore, the number of selected grids decreases for large $\lambda$ as before, providing us with more compressed result and easier understandings of the learning. 


\subsection{Logistic regression on binary images - a hard example}
\label{sec:logreg_im_difficult}
\begin{figure}[htbp]
  \centering
  \includegraphics[width=\hsize]{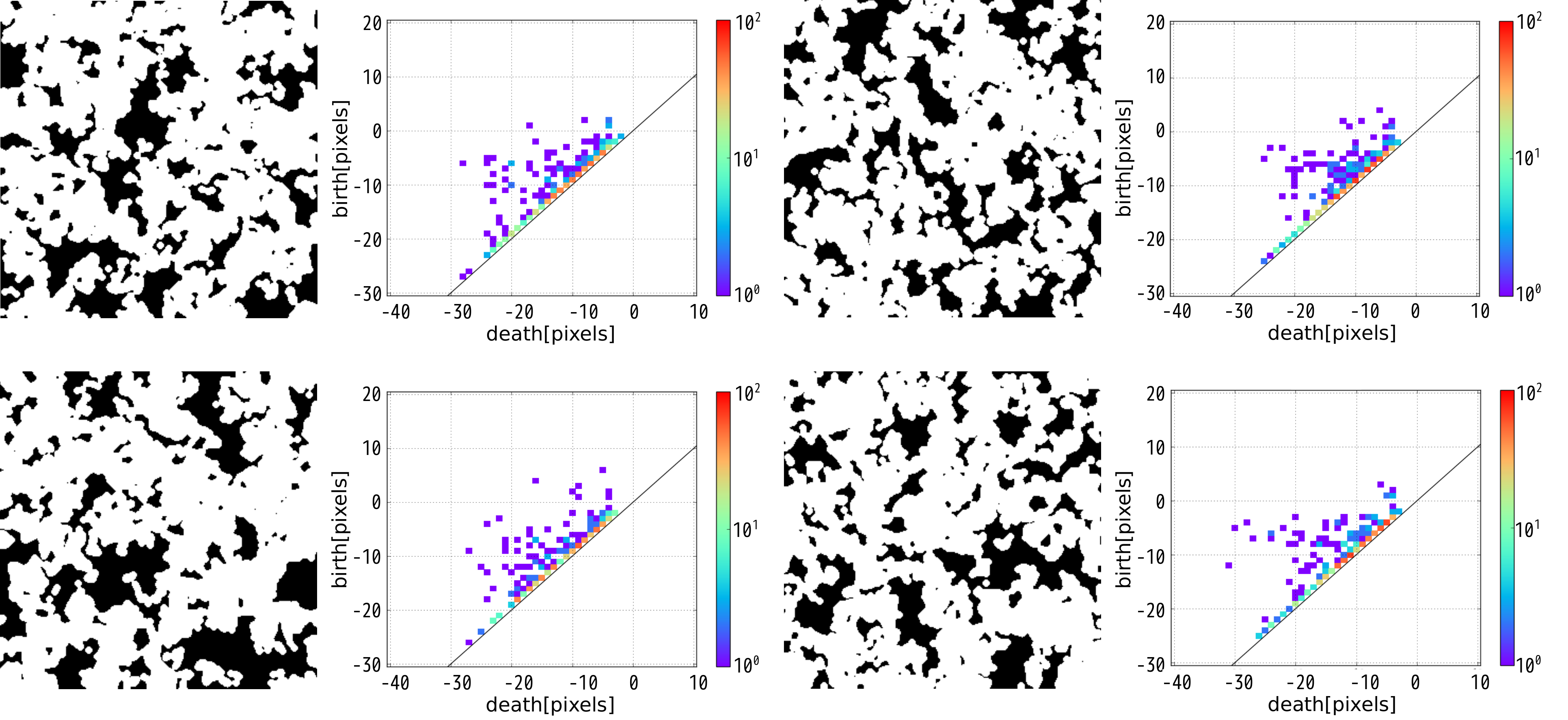}
  \caption{Input binary images and their 0th persistence diagrams. The left and right images are sampled from the parameter pairs (C) and (D), respectively. }
  \label{fig:input_im_B}
\end{figure}

Next, let us set the parameter pairs for generating random binary images so that the classification task becomes more difficult. Here, one parameter pair (C) is set to be $N=160, S=34$ and the other pair (D) is $N=270, S=18$. Figure~\ref{fig:input_im_B} shows the sample input data (left: (C), right: (D)).

In this example, it seems difficult to distinguish two parameter pairs based on our intuition. In fact, simple descriptors such as 
the average numbers of connected components and white pixels do not work at all in this case.

\begin{figure}[htbp]
  \centering
  \includegraphics[width=0.65\hsize]{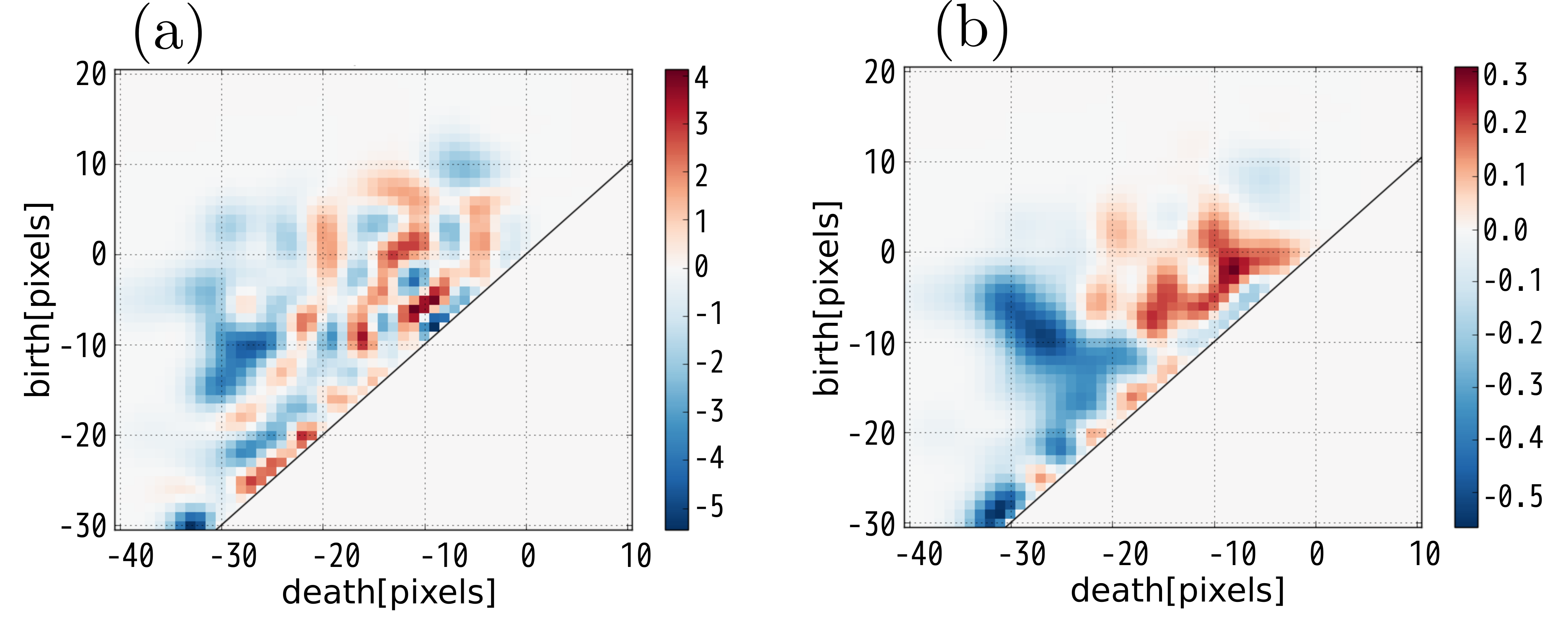}
  \caption{The reconstructed persistence diagrams using the $\ell^2$-regularization for the parameter pairs (C) and (D). (a) $\lambda = 0.0059948$ (chosen by cross validation), (b) $\lambda=1$.
  }
  \label{fig:output_im_B}
\end{figure}

The setting for the classification is the same as before, i.e., 2$\times$200 for training and 2$\times$100 for the test, and we assign 0 and 1 to (C) and (D), respectively. In this case, the score on the test set is 0.92 (baseline: 0.5). Figure~\ref{fig:output_im_B} (a) shows the reconstructed persistence diagram as the learned result using the $\ell^2$-regularization with the weight $\lambda = 0.0059948$ determined by the cross validation. 

In this learning, the distribution of the reconstructed persistence diagram looks complicated to observe clear features. Hence, let us increase  the weight parameter $\lambda$  for simplifying the distribution. 
Figure~\ref{fig:output_im_B} (b) shows the result with $\lambda=1$, where its score of the learning is 0.91. It should be noted that, although the score becomes only a little worse, the distribution turns out to be simple enough to conclude that the red area is dominant in the region with the birth scale $>-20$. From this simplification, we can explicitly obtain geometric reasonings for this classification in a similar way to Section \ref{sec:easy}.

Now we compare our method to other standard methods for image classifications. The list of methods and those scores are summarized in Table~\ref{tab:logreg_results}. These demonstrations show that persistence images with the logistic regression and support vector machine have better accuracy than the others. In particular, we note that the performance of our method is better than the bag of keypoints approach with sift feature, which is one of the standard techniques for image classifications \citep{sift,bow,bow-classification,bow-strategy}.
This is because such standard image classification techniques are developed mainly for clearly distinguishable and well-structured objects such as photos of faces, artificial objects, or landscapes, and not for images like this example. This suggests that our approach using persistence diagrams has an advantage to disordered images, which are frequently observed in materials science data \citep{iron}.

\begin{table}[htbp]
  \centering
  \begin{tabular}{|p{8cm}|l|} \hline
    Method & Mean accuracy \\ \hline\hline
    PI, logistic regression, $\ell^2$-penalty & 0.92 \\ \hline
    Bag of keypoints using sift with grid sampling,
    SVM classifier with $\chi^2$ kernel & 0.85 \\ \hline
    \# of connected components of black pixels & 0.73 \\ \hline
    \# of connected components of white pixels & 0.50 \\ \hline
    \# of white pixels & 0.50 \\ \hline
  \end{tabular}
  \caption{Performance comparison (PI: persistence image, SVM: support vector machine). 
  }
  \label{tab:logreg_results}
\end{table}


\subsection{Logistic regression on point clouds}
In this example, the input point clouds are prepared from two different random
point processes; one is Poisson point process (PPP) and
the other is Ginibre point process (GPP) on a unit disk.
It is known that PPP has no interaction between points, while GPP has a repulsive interaction. The parameters for these two point processes are
adjusted so that the mean number of points on the disk is 30. 
The task in this example is to identify PPP or GPP for test point clouds. 
To this task, we apply our method  to the 1st persistence diagrams with  the $\ell^2$-logistic regression.

\begin{figure}[htbp]
  \centering
  \includegraphics[width=0.8\hsize]{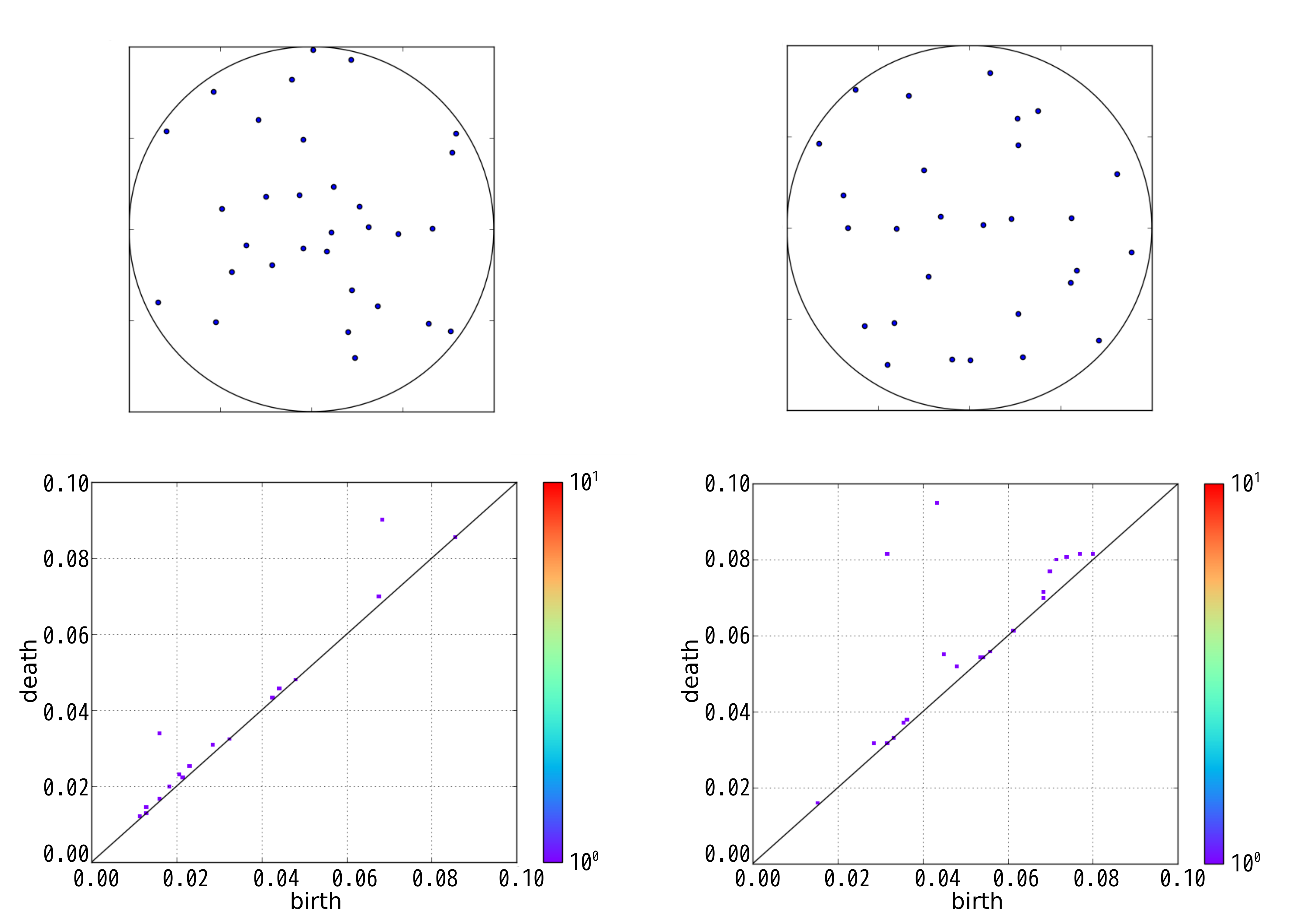}
  \caption{Random point clouds (Left: PPP, Right: GPP) and their 1st persistence diagrams}
  \label{fig:input_PPP_GPP}
\end{figure}

Figure~\ref{fig:input_PPP_GPP} shows point clouds generated by  
PPP and GPP. The parameters of the persistence images are set to be 
$\sigma=0.003, C=80, p = 1.0$ and the mesh for the discretized persistence diagrams is obtained by dividing the square $[0, 0.15]^2$ into $150\times 150$ grids.  
For each point process, 300 point clouds are generated (total $2 \times 300$)
and 200 of these point clouds are sampled as a training set (total $2 \times 200$), where we assign $0$ and $1$ for PPP and GPP, respectively. The remaining $2 \times 100$ point clouds are used as a test set for evaluation. 
Here, the weight parameter $\lambda$ of the $\ell^2$-regularization is determined by the cross validation. The score of the learned result  is 0.94.

\begin{figure}[htbp]
  \centering
  \includegraphics[width=0.8\hsize]{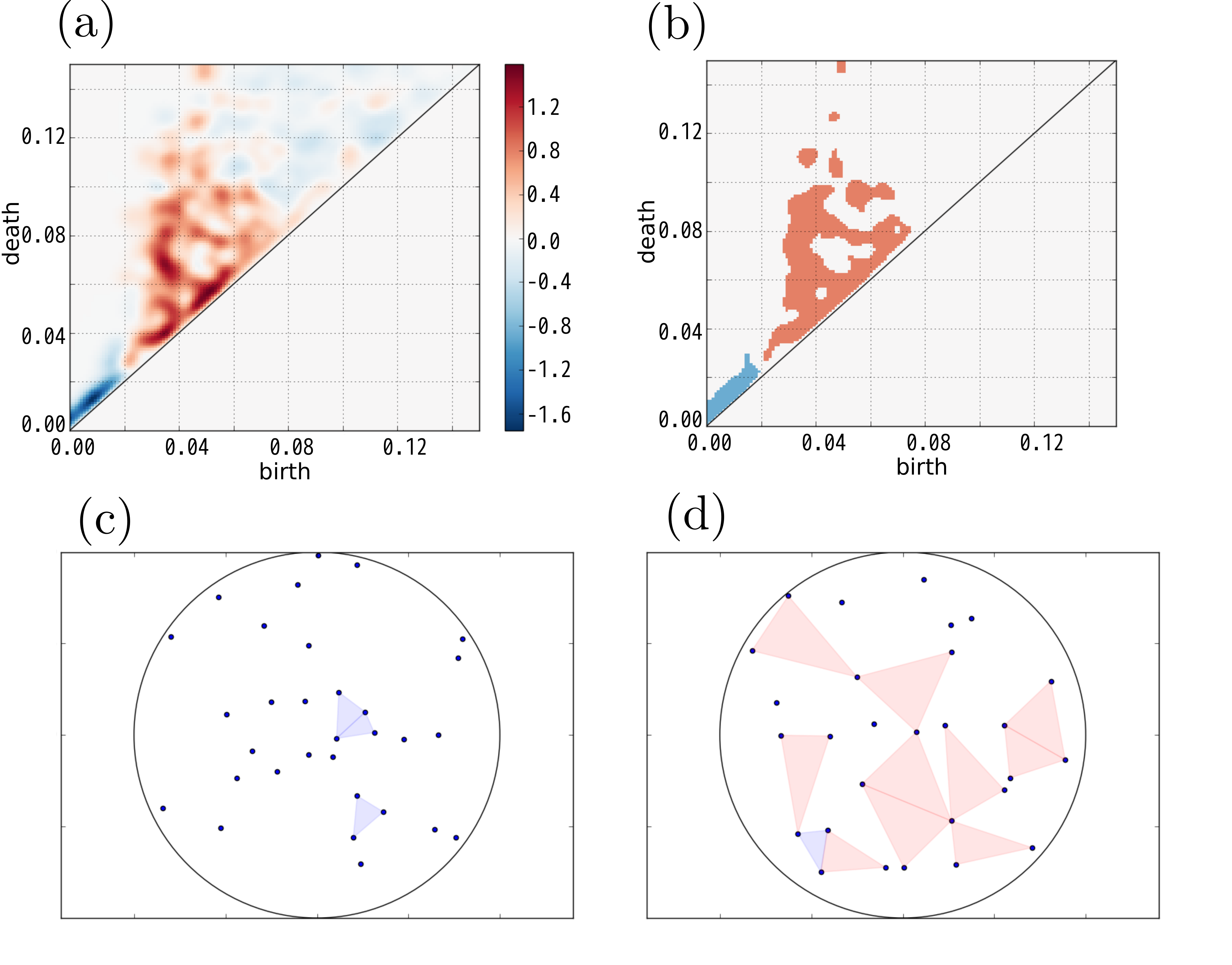}
  \caption{
    (a) The reconstructed persistence diagram. (b) A thresholding of (a).
    (c) The death positions (triangles) in PPP. (d) The death positions (triangles) in GPP.}
  \label{fig:output_pc_logreg}
\end{figure}

Figure~\ref{fig:output_pc_logreg} (a) shows the reconstructed persistence diagram from the learned vector $w$ and (b) shows the positive and negative areas of (a) with a certain threshold. 
Recall that, from the $0/1$ assignment,  the generators in the blue (resp. red) area contributes to classifying into PPP (resp. GPP).
From the learned persistence diagram, we observe that the red area  is located on the region with large birth values. This is consistent to the fact that GPP has a repulsive interaction, and hence it prevents the point cloud from constructing rings with small birth values. Figure~\ref{fig:output_pc_logreg} (c) and (d) show the death positions of the generators in the blue and red areas of (b) with the same colors, where (c) (resp. (d)) corresponds to PPP (resp. GPP). Similarly to the discussion in Figure 
\ref{fig:output_im_A_1}, these death positions express characteristic geometric features used for learnings more explicitly.

We remark that PPP and GPP can also be distinguished by using other descriptors such as average nearest neighbor distances. An advantage of our method is that we do not need any prior knowledge, providing us with more universal method compared to problem-specific descriptors. In fact, the analysis using average nearest neighbor distance can be realized by the $0$th persistence diagram. 

\begin{figure}[htbp]
  \centering
  \includegraphics[width=0.8\hsize]{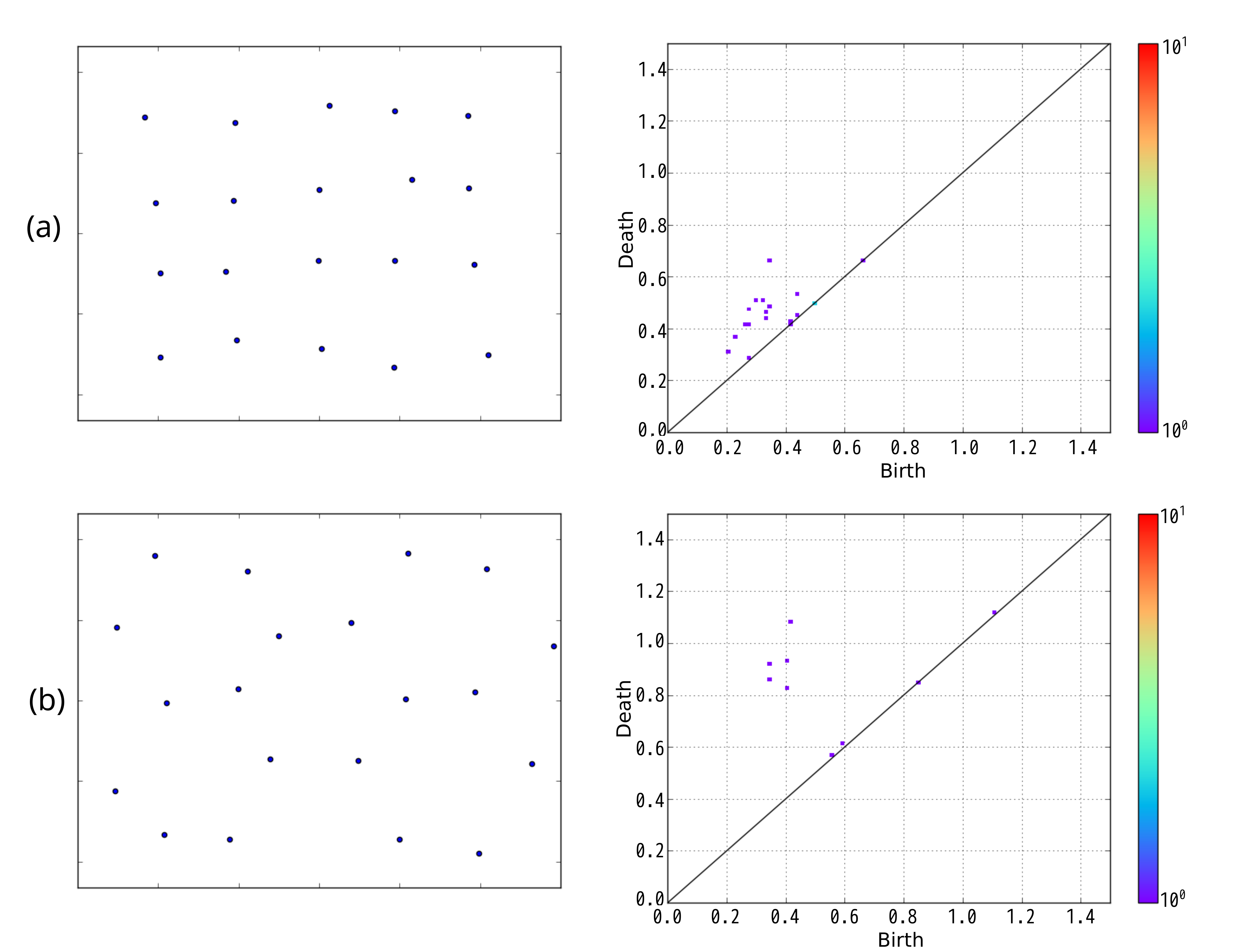}
  \caption{Input point clouds and their 1st persistence diagrams.
  (a) A square lattice with noise (b) A regular hexagonal lattice with noise.}
  \label{fig:pc-crystal-input}
\end{figure}

Now let us test another example for point clouds. The task is classifying two types of point clouds; one is a square lattice with Gaussian noise, and the other is a regular hexagonal lattice with Gaussian noise. Figure~\ref{fig:pc-crystal-input} shows
the input point clouds and their 1st persistence diagrams. In this example, the distance between two nearest neighbors is one for both cases, and hence it is difficult to distinguish these two types of point clouds using average nearest neighbor distances.

\begin{figure}[htbp]
  \centering
  \includegraphics[width=0.45\hsize]{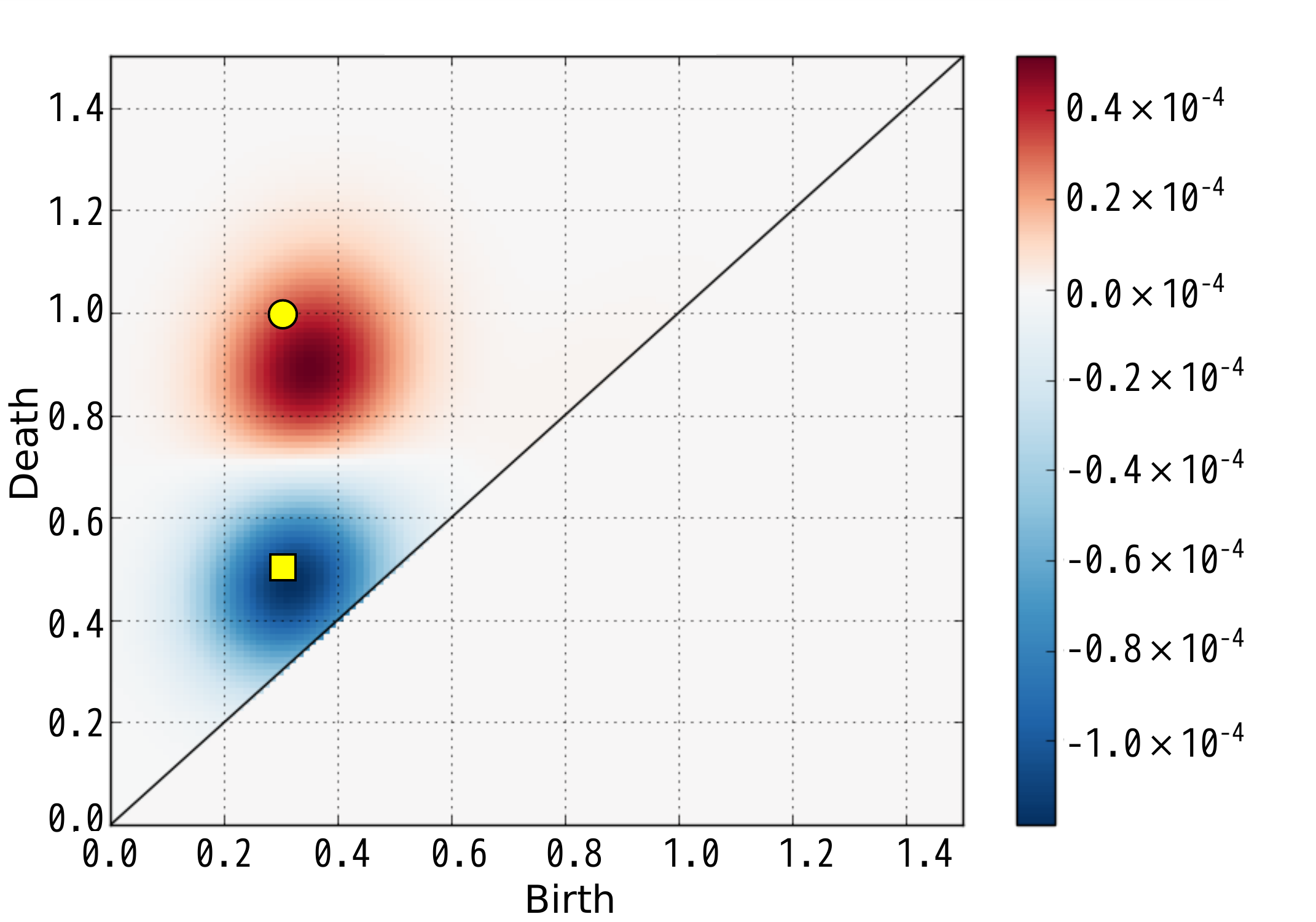}
  \caption{The reconstructed persistence diagram. The yellow circle (resp. rectangle)
  shows the birth-death pair corresponding to the regular hexagon (resp. square).}
  \label{fig:pc-crystal-output}
\end{figure}

We set the number of points to be 20 and the standard deviation of the noise to be $\sigma = 0.1$. 
Figure~\ref{fig:pc-crystal-output} shows the reconstruct persistence diagram from the learned vector $w$.  The yellow circle (resp. rectangle) in the diagram shows the birth-death pair of the regular hexagon (resp. square).  One interesting feature in this result is that the positive peak position in the reconstructed diagram is shifted to the diagonal from yellow circle. Probably this is because such a regular shape is optimal in order to leave from the diagonal, and many birth-death pairs in noisy hexagonal lattices tend to move toward the diagonal.

\subsection{Linear regression on binary images}
\begin{figure}[htbp]
  \centering
  \includegraphics[width=0.8\hsize]{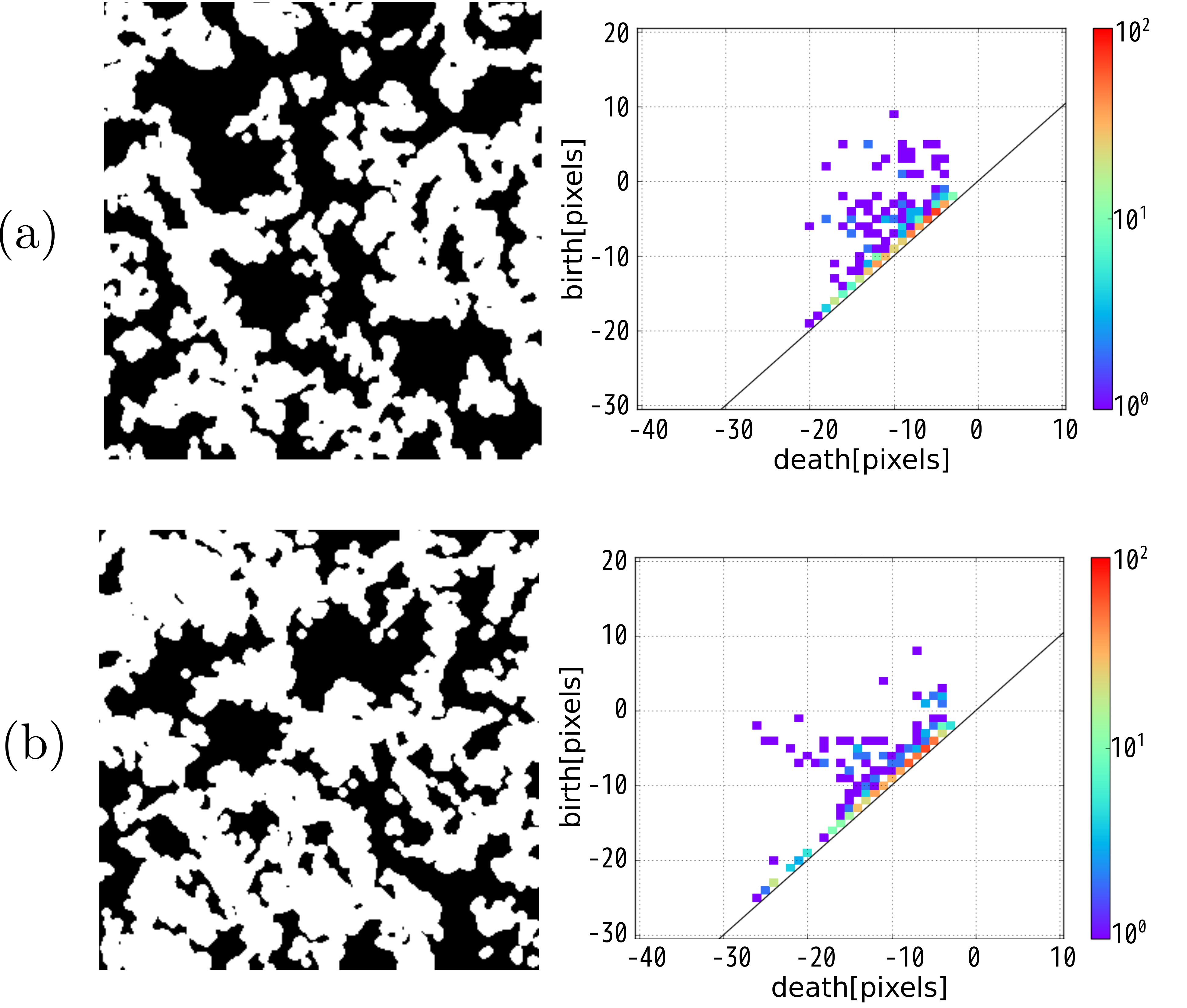}
  \caption{Sample images for the linear regression and
    their 0th persistence diagrams. (a) $S=21$, (b) $S=28$.}
  \label{fig:im_linreg_input}
\end{figure}

In this example, we examine the linear regression on binary images.
The input binary images are generated by Algorithm~\ref{alg:randomimage}
with $N=150$ and $S$ is randomly chosen from $\{20, 21, \cdots, 29\}$ uniformly.
The task is to determine the random parameter $S$ from images.
Figure~\ref{fig:im_linreg_input} shows sample images with $S=21$ and $S=28$ and those persistence diagrams. 

From the construction of Algorithm~\ref{alg:randomimage}, we know that $S$ controls the area of white pixels. Hence, to our task, we study the following descriptors:
\begin{enumerate}[(i)]
\item persistence image
\item the number of white pixels
\item the combination of (i) and (ii)
\end{enumerate}
as explanatory variables and compare these performances.
Here, the third descriptor means that the response variable $S$ is explained by the following model
\begin{align}
  S = v\cdot\textrm{(\# of white pixels)} + w \cdot \textrm{(PI)} + b +
      \textrm{(noise),} \label{eq:im_linearmodel} 
\end{align}
where $v, b \in \R$ and $w \in \R^n$ are unknown parameters and determined from a training set. For (i) and (iii), we apply both $\ell^2$- and $\ell^1$-regularizations. The weight parameter $\lambda$ of the regularization is determined by the cross validation.

The training set and test set consist of randomly generated 500 images and 100 images, respectively. The learned results are assessed using the $R^2$ coefficients of determination\citep{statistics} on the test set, which are shown in Table~\ref{tab:linreg_r2}.
As we observe, our methods (i) using $\ell^1$- and $\ell^2$-regularizations attain almost the same performance as (ii), while the combination (iii) improves the performance better. 

\begin{table}[htbp]
  \centering
  \begin{tabular}{|c|l|l|} \hline
    &Method & $R^2$ coefficient \\ \hline\hline
    \multirow{2}{*}{(i)} & PI with ridge ($\ell^2$) & 0.86 \\ \cline{2-3} 
    &PI with lasso ($\ell^1$) & 0.86 \\ \hline 
    (ii)&\# of white pixels & 0.88 \\ \hline 
    \multirow{2}{*}{(iii)}&Both with ridge & 0.93 \\ \cline{2-3} 
    &Both with lasso & 0.94 \\ \hline 
  \end{tabular}
  \caption{$R^2$ coefficients on the test set of the linear regression problem. These values become larger when the learned model gives better predictions.}
  \label{tab:linreg_r2}
\end{table}

\begin{figure}[htbp]
  \centering
  \includegraphics[width=0.8\hsize]{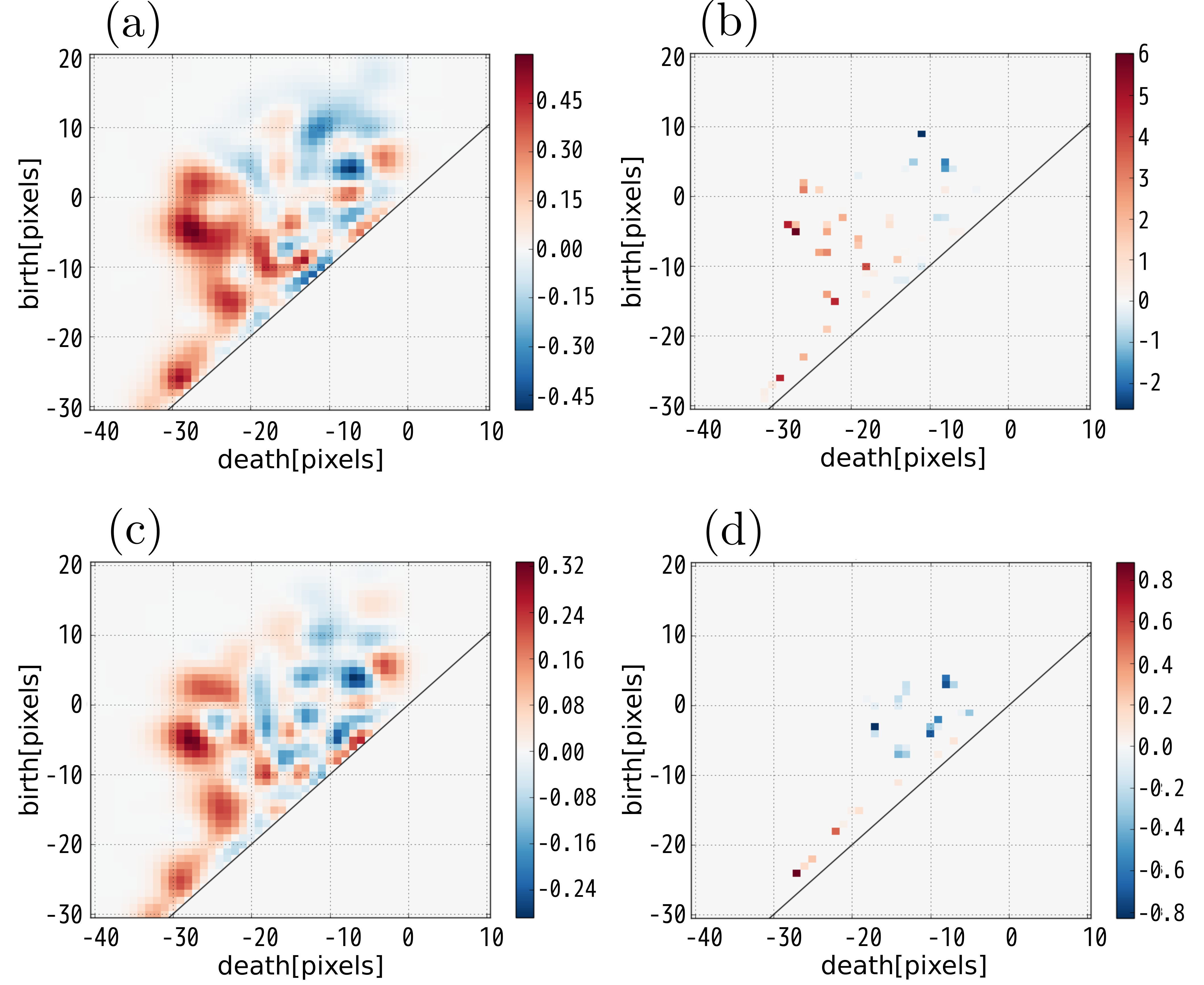}
  \caption{The reconstructed persistence diagrams. (a) PI with ridge. (b) PI with lasso. (c) Both with ridge. (d) Both with lasso.}
  \label{fig:ridge_lasso}
\end{figure}

Figures~\ref{fig:ridge_lasso} shows the reconstructed persistence diagrams obtained from (i) and (iii). By construction of our regression model, the areas with  positive (resp. negative) values on the diagrams positively (resp. negatively) contribute to the response variable $S$. Even in the linear regression model, we can observe the sparseness property for the $\ell^1$-regularization, which is useful for extracting the most essential features for the response variable $S$ from sample data.

From the mixed model (\ref{eq:im_linearmodel}), we can estimate the contributions of (i) and (ii) in (iii) for predictions. For example, the following prediction results applied to Figure ~\ref{fig:im_linreg_input} (a) and (b) with the $\ell^1$-regularization imply that the prediction mainly consists of
the term $v\cdot\textrm{(\# of white pixels)}$ and is modified negatively by
the term $w \cdot \textrm{(PI)}$. 
\begin{align*}
  \begin{array}{lllll}
    S & \approx
        \mbox{(prediction of $S$)} & = v\cdot\textrm{(\# of white pixels)} &+ w \cdot \textrm{(PI)} &+  b  \\
    21 &\approx 20.628 &= 30.272 &+  (-5.917) &+  (-3.728) \\
    28 &\approx 27.959 &= 35.718 &+  (-4.031) &+  (-3.728)   
  \end{array}
\end{align*}

\begin{figure}[htbp]
  \centering
  \includegraphics[width=1.0\hsize]{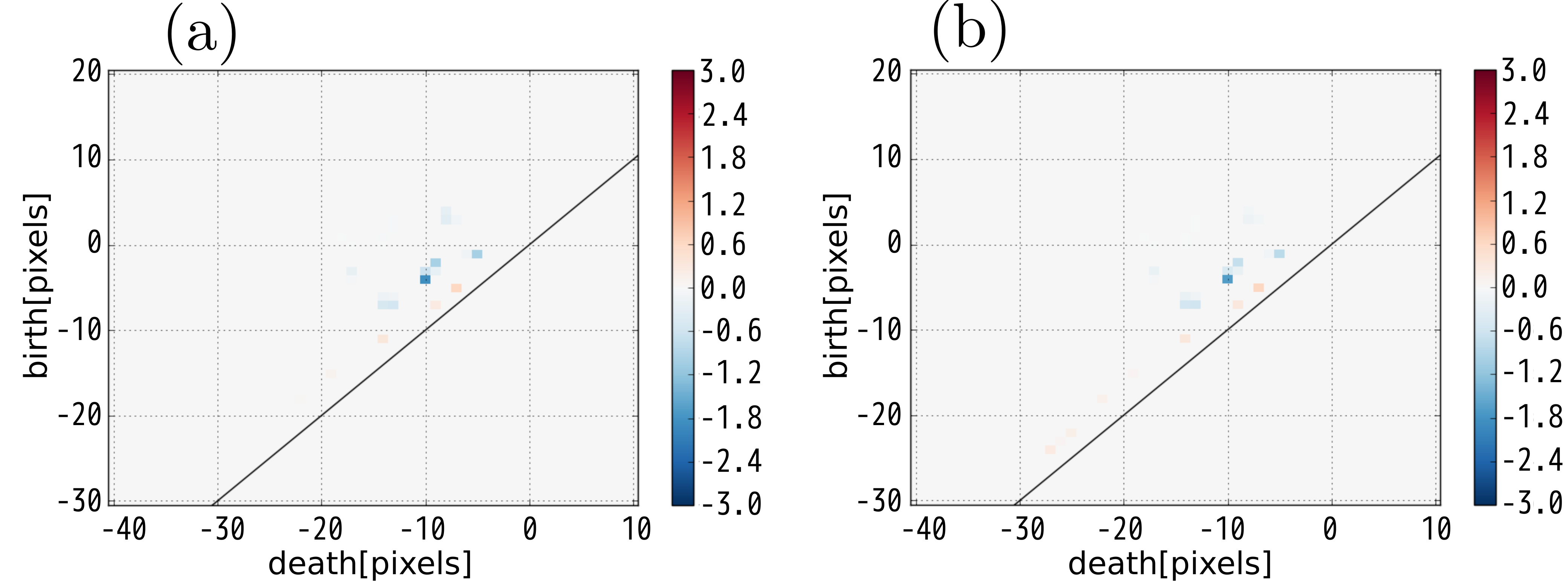}
  \caption{The weighted persistence diagrams for
    Figure~\ref{fig:im_linreg_input} (a) and (b).}
  \label{fig:im_lasso_weights}
\end{figure}

Furthermore, by showing the weighted persistence diagram $(w_ix_i)_{i=1}^{n}$ for the test persistence diagram $x$, we can explicitly clarify the important generators for modifications. 
Figures~\ref{fig:im_lasso_weights}  shows the weighted persistence diagrams of Figure~\ref{fig:im_linreg_input} (a) and (b), and in this case, we find that generators around $(-10, -4)$ effectively work for predictions of $S$.

For applications in materials science, $S$ can be regarded as a certain physical quantity such as conductivity of battery materials. Then, by this approach, we can identify  geometric structures in the images which most effectively affect that physical quantity. 

\section{Conclusion}\label{sec:conclusion}
In this paper, we have proposed a unified method by combining persistence images and linear machine learning models with the ability to study the inverse problem in the original data space. One of the important properties of our method is that a persistence diagram is obtained as a learned result. From such a reconstructed persistence diagram and the inverse analysis using birth/death positions, we can explicitly characterize significant  geometric features embedded in dataset. We have also presented sparse persistence diagrams as an important concept of machine learnings in topological data analysis. 

Although we applied our method to linear regressions and logistic regressions, it is obviously not limited to them, and many other linear machine learning models such as support vector machine with a linear kernel and elastic nets are also applicable. Moreover, we can similarly apply our method to point clouds and cubical sets in higher dimensions.

The proposed method is recently applied to several practical problems. For example, in the paper \citep{iron}, the authors develop  a method for predicting locations of micro cracks generated by reduction reaction process of iron ore sinters. In that application, they apply the persistence images with the $\ell^1$-linear regression to a huge amount of X-CT images, and select the crack areas as a response variable. Then, it follows that the reconstructed persistence diagram from the learned vector identifies  generators which have significant effects on crack formations, and hence, by studying their birth/death positions, we can explicitly detect the location of micro cracks. We believe that the same analysis is also useful to other problems dealing with large amount of images such as pathology.

\appendix

\section{Algorithm for generating random images}
\label{sec:randomimage}

The algorithm for generating random binary images is given 
by Algorithm~\ref{alg:randomimage}. It consists of six parameters, $W, N, S\in \mathbb{N}, \sigma_1 > 0, \sigma_2 > 0$, and $t >0$.
The area of white pixels in the generated image
is given by the orbits of the Brownian motion of $N$ particles on a flat torus with the size $W \times W$. The parameters $S$ and $\sigma_1$ determine the length of each orbit and $\sigma_2$ and $t$ determine the radii of particles. In this paper we fix $W=300$, $\sigma_1 = 4$, $\sigma_2 = 2$, $t = 0.01$, and only $N$ and $S$ are changed. When $N$ and $S$ become larger, the generated image tend to have more white pixels. 

These kinds of random images are frequently obtained by experimental measurements in materials science such as X-CT and TEM \citep{iron}. These seemingly disordered images are supposed to be utilized for materials informatics, and one of the motivations of this paper is to develop a universal framework for this purpose. 

\begin{algorithm}[h!]
  \caption{Generate a random binary image}
  \label{alg:randomimage}
    \begin{algorithmic}
    \Procedure{Gen-Image}{$W, N, S, \sigma_1, \sigma_2, t$}
    \State Let $T$ be $[0, W]\times [0, W]$
    \For{$n=1,\ldots,N$}
      \State Take $x_{n,1}$ uniformly randomly on $T$
      \For{$s=1,\ldots,S$}
        \State Take $d_1$ and $d_2$ randomly from ${\mathcal N}(0, \sigma_1)$
        \State $x_{n,s+1} \gets x_{n,s} + (d_1, d_2) \mod W \times W$
      \EndFor
    \EndFor
    \State $H$ $\gets$ The Histogram of $\{x_{n,s}\}$ with $W\times W$
      mesh on $T$
    \State Apply Gaussian filter to $H$ with the standard deviation $\sigma_2$ and set the result to $\tilde{H}$
    \State \Return The Binary image from $\tilde{H}$ by thresholding with $t$
    \EndProcedure
  \end{algorithmic}

\end{algorithm}


\end{document}